\documentclass[12pt]{article}
\usepackage{longtable,supertabular}
\usepackage{amsmath}
\usepackage{amssymb}
\usepackage{latexsym}
\usepackage{a4wide,color}
\usepackage{multirow}
\usepackage{ntheorem}
\usepackage[square,numbers,sort]{natbib}
\usepackage[title]{appendix}
\usepackage{array}
\usepackage{makecell}

\title{\bf A Simple Explicit Construction of Smaller Strong Blocking Sets}
\author{Hao Chen \thanks{Hao Chen is with the College of Information Science and Technology, Jinan University, Guangzhou, Guangdong Province, 510632, China (e-mail: haochen@jnu.edu.cn). This research was supported by NSFC Grant 62032009.}}

\begin{document}
	
	\maketitle
	\begin{abstract}
		Alon, Bishnoi, Das and Neri constructed sequences of strong blocking sets in $(k-1)$-dimension projective space $P({\bf F}_q, k-1)$ with sizes at most $Cqk$, where $k$ goes to the infinity and $C$ is a fixed positive constant not depending on $q$, in their paper published in Transactions AMS, vol. 277, pp. 5389-5410, 2024. This is equivalent to a sequence of minimal $[n_i,k_i]_q$ codes with positive rate $\frac{k_i}{n_i} \geq c$, where $k_i$ goes to the infinity and $c$ is a fixed positive constant only depending on $q$. They solved one of the main open problems on minimal codes. The Alon-Bishnoi-Das-Neri construction was based on asymptotically good sequences of algebraic geometric codes and expander graphs. However, the constant $C$ is large for small prime power $q$. \\

In this paper, based on a simple approach using $(q-1)$-power intersecting algebraic geometric codes via concatenation to construct minimal codes, we construct explicit strong blocking sets in $(k-1)$-dimension projective space with sizes around $2(q^2-1)k$, equivalently, explicit minimal $[n,k]_q$ codes with length around $2(q^2-1)k$. This is an essential improvement on the previous result of Alon-Bishnoi-Das-Neri, when $q \leq 11$. For large $q$, strong blocking sets and minimal codes in our construction are worse than the previous construction of Alon-Bishnoi-Das-Neri. However, our approach is much simpler and does not depend on expander graphs.\\

		 Alon-Bishnoi-Das-Neri construction.
	\end{abstract}
	
	\newpage

\section{Introduction}

\subsection{Strong blocking sets}

Let $P({\bf F}_q, k-1)$ be the $(k-1)$-dimensional projective space over the finite field ${\bf F}_q$. Blocking sets in projective spaces have been studied for many years, see \cite{BS,Beu,BB}. A subset $S$ in the $(k-1)$ dimension projective space $P(q,k-1)$ over ${\bf F}_q$ is called a $t$-fold $s$-blocking set, if for each subspace $H \subset P(q,k-1)$ of codimension $s$, the intersection $S \bigcap H$ has at least $t$ points. A subsect $S \subset P({\bf F}_q, k-1)$ is called cutting $s$-blocking set, if for each subspace $H \subset P(q,k-1)$ of codimension $s$, the intersection $S \bigcap H$ spans the whole space $H$. The cutting $1$-blocking set is also called strong blocking set. It is obvious that a cutting $s$-blocking set is a $(k-s)$-fold $s$-blocking set. From the classical result by Beutelpacher in \cite{Beu}, if $t \leq q$, then the size of a $t$-fold $s$-blocking set $S$ is at least $$|S| \geq t(q^s+q^{s-1}+\cdots+1).$$ A new lower bound on sizes of $t$-fold $s$-blocking sets without assuming $t \leq q$, which is stronger than Beutelpacher's result,  was given in our recent paper \cite{Pan1}. For strong block sets, this lower bound was strengthened recently to a lower bound $$|S| \geq (q+1)(k-1)$$ in \cite{Alfarano}. On the other hand, from a probabilistic argument, it was proved that there exists a strong blocking set of the size at most $$|S| \leq \left\lceil \frac{2}{1+\frac{1}{(q+1)^2ln q}}(k-1)\right\rceil (q+1)$$ for $q \geq 3$, and $$|S| \leq \frac{2k-1}{log_2\frac{4}{3}}$$ for $q=2$, see \cite{Mik,Alfarano3,Bishnoi,Hegar}. In the case of $char({\bf F}_q)>k$, strong blocking sets in $P({\bf F}_q, k-1)$ of sizes $O(qk)$ was constructed in \cite{FS}. For the construction of strong blocking sets in $P({\bf F}_q,k-1)$ of sizes $O(k^2)$, we refer to \cite{Davy,Alfarano}.\\

\subsection{Error-correcting codes}

The Hamming weight $wt({\bf a})$ of a vector ${\bf a}=(a_0, \ldots, a_{n-1}) \in {\bf F}_q^n$ is the cardinality of its support, $$supp({\bf a})=\{i: a_i \neq 0\}.$$ The Hamming distance $d({\bf a}, {\bf b})$ between two vectors ${\bf a}$ and ${\bf b}$ is $d({\bf a}, {\bf b})=wt({\bf a}-{\bf b})$. Then ${\bf F}_q^n$ is a finite Hamming metric space. The minimum (Hamming) distance of a code ${\bf C} \subset {\bf F}_q^n$ is, $$d({\bf C})=\min_{{\bf a} \neq {\bf b}} \{d({\bf a}, {\bf b}),  {\bf a} \in {\bf C}, {\bf b} \in {\bf C} \}.$$ A linear code of dimension $k$ in ${\bf F}_q^n$ is denoted by a linear $[n,k]_q$ code. There are some coding upper bounds on parameters, such as Singleton bound $$d \leq n-k+1$$ for a linear $[n,k,d]_q$ codes. One of the main tasks in the theory of error-correcting codes is to construct error-correcting codes with large $k$ and $d$, when $n$ and $q$ are given.\\

Let ${\bf C} \subset {\bf F}_q^m$ be a linear code over ${\bf F}_{q^m}$, We denote the linear code over ${\bf F}_{q^m}$, spanned by all vectors $${\bf c}^{q-1}=(c_0^{q-1}, \ldots, c_{n-1}^{q-1}),$$ where ${\bf c}=(c_0,\ldots, c_{n-1})\in {\bf C}$ takes all nonzero codewords of ${\bf C}$, over ${\bf F}_{q^m}$ by ${\bf C}^{q-1}$. As Schur products of linear codes, when ${\bf C}$ is an evaluation code, then it is not hard to get a linear code containing the $(q-1)$-power code of ${\bf C}$.\\

\subsection{Minimal codes and intersecting codes}

Minimal binary linear codes were studied in \cite{CL} in the name linear interesting codes, also see \cite{Borello}. In \cite{AB}, minimal codewords and minimal codes were introduced for general $q$-ary linear codes and studied systematically.  A nonzero codeword ${\bf c}$ of a linear $[n,k]_q$ code ${\bf C}$ is called minimal, if for any codeword ${\bf c}'$ satisfying $supp({\bf c}') \subseteq supp({\bf c})$, then there is a nonzero $\lambda \in {\bf F}_q$, such that, ${\bf c}'=\lambda {\bf c}$. A linear code ${\bf C}$ is called minimal if each nonzero codeword is minimal. In \cite{AB}, the Ashikhmin-Barg criterion $\frac{w_{min}}{w_{max}} > \frac{q-1}{q}$ was proposed as a sufficient condition for a $q$-ary linear code to be minimal, where $w_{min}$ and $w_{max}$ are minimum and maximum weights of the code. \\

A linear code ${\bf C} \subset {\bf F}_q^n$ is called intersecting, if the supports of any two nonzero codewords have no empty intersection. That is, $$|supp({\bf c}) \bigcap supp({\bf c}')|\geq 1,$$ for any two nonzero codewords of ${\bf C}$. For binary linear codes, intersecting codes are equivalent to minimal codes. When $q \geq 3$, minimal $q$-ary linear codes are intersecting. However, there are $q$-ary non-minimal intersecting codes, see \cite{CL,CZ,Borello}. A linear code over ${\bf F}_{q^m}$ is called $(q-1)$-power intersecting if ${\bf C}^{q-1}$ is an intersecting $q^m$-ary code.\\

\subsection{Strong blocking sets and minimal codes}

The following geometric characterization of minimal linear codes was given in \cite{Tang,Alfarano1}, also see \cite{Hegar,Alon}. Let ${\bf G}=({\bf g}_1,\ldots, {\bf g}_n)$ be the generator matrix of a dimension $k$ linear code ${\bf C} \subset {\bf F}_q^n$, where ${\bf g}_1, \ldots, {\bf g}_n$ are nonzero vectors in ${\bf F}_q^k$. Let $S$ be the subset of $n$ points in $P(q,k-1)$ corresponding to ${\bf g}_1, \ldots, {\bf g}_n$. Then ${\bf C}$ is a minimal code if and only if $S$ is a strong blocking set. In \cite{Alfarano2}, this geometric characterization of minimal codes was generalized to minimal codes in the rank matric.\\

From the geometric characterization of minimal codes, the lower bounds on sizes on strong blocking sets can be translated to lower bounds on lengths of minimal linear codes. The minimal length $m(k,q)$ of a linear $q$-ary dimension $k$ minimal code has to satisfy $$m(k,q) \geq (k-1)(q+1).$$ Based on a probabilistic argument, $$m(k,q) \leq \left\lceil \frac{2}{1+\frac{1}{(q+1)^2ln q}}(k-1)\right\rceil (q+1)$$ for $q \geq 3$, and $$m(k,2) \leq \frac{2k-1}{log_2\frac{4}{3}}$$ for $q=2$ was proved, we refer to \cite{Mik,Alfarano3,Bishnoi,Hegar} for such kind of upper bound on $m(k,q)$. It was proved in \cite{Bishnoi} that $$m(k,q) \leq (q+1) \cdot \frac{2k}{log_q\frac{q^4}{q^3-q+1}}.$$ However, their proof is not constructive, therefore, it is interesting and important to construct explicit small strong blocking sets and short minimal codes. We refer to \cite{Scotti} for improvements on lower bounds on $m(k,q)$. \\

A ternary code ${\bf C} \subset {\bf F}_3^n$ is called trifferent, if and only if, for any three distinct codewords ${\bf x}, {\bf y} $ and ${\bf z}$ of the code ${\bf C}$, there exists an coordinate position $1\leq i \leq n$ such that, $\{x_i,y_i,z_i\}={\bf F}_3$. It was proved in \cite{Bishnoi} that a linear ternary code is trifferent, if and only if it is a minimal ternary code. For a given length $n$, there are previous results on the size of (linear) trifferent code. In \cite[Theorem 1.7]{Bishnoi}, it was proved that sequences of linear trifferent codes of the rate at least $\frac{23}{312}$ can be constructed explicitly. Their proof in \cite{Bishnoi} depends on an argument of Alon in \cite{Alon1}.\\

\subsection{Recent advances}

In \cite[Corollary 3.3]{Bartoli}, by the using of outer codes satisfying the Ashikhmin-Barg condition, linear minimal $q$-ary codes with the rate at least $\frac{4}{q^4}$, were constructed. Or equivalently, strong blocking sets in $P({\bf F}_q, k-1)$ with the size at most $\frac{q^4k}{4}$, were constructed. The recent breakthrough in \cite{Alon} gave constructions of strong blocking sets in $P({\bf F}_q,k-1)$ with sizes at most $Cqk$, and $k$ dimension minimal codes of length at most $Cqk$, by the using or asymptotically good algebraic geometric codes and expander graphs. Their constant $C$ is universal and does not depend on $q$. As cited from \cite{Alon}, "The main problem is to construct minimal codes of a given dimension $k$ of the shortest possible length $n$."  Notice that the minimum weight $d$ of a minimal linear $[n,k]_q$ code is at least $(q-1)(k-1)+1$, see \cite{Alfarano}. Then the constructed sequences of positive rate minimal codes in \cite{Alon} are automatically asymptotically good sequence of $q$-ary codes. As commented in the Abstract of \cite{Alon}, "This solves one of main open problems on minimal codes".\\

However, their positive constant $C$ is large, when $q=2$, it was proved that $C$ is around $3 \times 116$, for $q=3$, their $C$ is around $4 \times 77$, and for general prime power $q$, their $C$ is around $20(q+1)$, see
\cite[Table 3]{Alon}. In the case $q=2$, this positive constant was improved in \cite{Borello} to around $5.83$, since binary minimal
 codes are equivalent to binary intersecting codes, see \cite[Table 3]{Borello}. For general prime power $q \geq 3$, results in \cite[Table 3]{Borello} cannot be applied to minimal codes, since intersecting codes are not minimal codes when $q \geq 3$. When $q=2$, and $k$ is small, some short minimal codes were constructed in \cite{Kurz,CZ}. For the case of $q=3$, ternary minimal codes are equivalent to trifferent codes, we refer to \cite{Bishnoi} for short ternary minimal codes. We refer to \cite{Scotti1} for a nice survey on recent advance on strong blocking sets and minimal codes.\\

A ternary code ${\bf C} \subset {\bf F}_3^n$ is called trifferent, if and only if, for any three distinct codewords ${\bf x}, {\bf y} $ and ${\bf z}$ of the code ${\bf C}$, there exists an coordinate position $1\leq i \leq n$ such that, $\{x_i,y_i,z_i\}={\bf F}_3$. It was proved in \cite{Bishnoi} that a linear ternary code is trifferent, if and only if it is a minimal ternary code. For a given length $n$, there are previous results on the size of (linear) trifferent code. In \cite[Theorem 1.7]{Bishnoi}, it was proved that sequences of linear trifferent codes of the rate at least $\frac{23}{312}$ can be constructed explicitly. Their proof in \cite{Bishnoi} depends on an argument of Alon in \cite{Alon1}.\\

\subsection{Main results}

The main purpose of this paper is to construct smaller strong block sets than \cite{Alon} for small prime power $q \leq 11$, or equivalently, shorter minimal $q$-ary linear codes with better rates. Our new construction giving new longer minimal codes over ${\bf F}_q$ via shorter minimal codes over ${\bf F}_q$ and $(q-1)$-power intersecting codes over ${\bf F}_{q^m}$. Notice that the construction is not restricted to the case $q=2$ as in \cite{Borello}. For general prime power $q$, the constant $C$ is improved to around $$\frac{n(q^m-1)(q-1)}{m(q^m-3)}\approx 2(q^2-1)$$ if there exists a minimal $[n,2m]_q$ minimal code, where $m$ can be arbitrary positive integer $m \geq 2$. This is worse than that in \cite{Alon} for large prime power $q$. When $q$ is a small prime power satisfying $q\leq 11$, our constant is better than that of \cite{Alon}. For example, when $q=3$, there is a minimal $[15,4]_3$ code given in \cite[Page11]{Kurz}, then the constant $C$ is improved to $\frac{15 \cdot 8 \cdot 2}{2 \cdot 5}=20 << 4 \times 77$, which is much better than the result in \cite[Table 3]{Alon}. For general large prime power, our approach is much simpler and needs not expander graphs.\\

In this paper, we prove the following results.\\

{\bf Theorem 1.1.} {\em Let ${\bf C} \subset {\bf F}_{q^m}^N$ be a $q^m$-ary $K$-dimensional linear code such that ${\bf C}^{q-1}\subset {\bf F}_{q^m}^N$ is an intersecting code. We assume that there is an explicit $[n,m]_q$ minimal code in ${\bf F}_q^n$. Then we can construct a minimal linear $[Nn,Km]_q$ code in ${\bf F}_q^{Nn}$ explicitly.}\\

{\bf Theorem 1.2} {\em Suppose that $q$ is an arbitrary prime power and $m \geq 2$ is a positive integer. We assume that there exists an explicit $[n,2m]_q$ minimal code. Then we construct a sequence of strong blocking sets $S_i \subset P({\bf F}_q, k_i-1)$, where $k_i$ goes to the infinity, such that $$\lim\frac{|S_i|}{k_i}\leq \frac{n(q^m-1)(q-1)}{m(q^m-3)}.$$ Equivalently, we construct sequences of minimal $q$-ary codes with rate at least $\frac{m(q^m-3)}{n(q-1)(q^m-1)}$ explcitly.}\\

{\bf Corollary 1.1.} {\em For $q=2$, we construct a sequence of linear binary minimal codes with the rate $\frac{145}{961} \approx \frac{1}{6.6276}$ explicitly.}\\

This is a little worse than $\frac{1}{5.8334}$ in \cite[Table 3]{Borello}, since binary intersecting codes are equivalent to binary minimal codes.\\

{\bf Corollary 1.2.} {\em For $q=3$, we construct a sequence of linear trifferent (ternary minimal) codes with the rate $\frac{1}{20}$ explicitly.}\\

This is a little worse than $\frac{23}{312}\approx \frac{1}{13.56}$ in \cite[Theorem 1.7]{Bishnoi}, since an argument in \cite{Alon1} was used in the proof of Theorem 1.7 of \cite{Bishnoi}.\\

However, for $q \geq 4$, there is no better explicit construction of minimal $q$-ary codes, than the results given in \cite[Table 3]{Alon}. We have the following results on explicit construction of small size strong blocking sets or short minimal codes for small prime powers $q \geq 4$.\\

{\bf Corollary 1.3.} {\em For $q=4$, we construct a sequence of linear $4$-ary minimal codes with the rate at least $\frac{13}{405}$ explicitly.}\\

{\bf Corollary 1.4.} {\em For $q=5$, we construct a sequence of linear $5$-ary minimal codes with the rate at least $\frac{11}{864}$ explicitly.}\\

This paper is organized as follows. In Section 2, we construct algebraic geometric $q^m$-ary codes ${\bf C}$, such that ${\bf C}^{q-1}$ is intersecting.  In Section 3, we give a simple construction of strong blocking sets and minimal codes from $(q-1)$-power intersecting algebraic geometric codes, via the concatenation. In Section 4, we calculate the asymptotic upper bound of $\lim_{k \longrightarrow \infty}\frac{m(k,q)}{k}$. \\

\section{$(q-1)$-power intersecting algebraic geometric codes}

We recall some basic facts about algebraic geometric codes, see \cite{TV}. Let ${\bf X}$ be an projective smooth algebraic curve of genus $g$ defined over the finite field ${\bf F}_q$. Let ${\bf P}=\{P_1,\ldots, P_n\}$ be a set of ${\bf F}_q$-rational points of ${\bf X}$ and ${\bf G}$ be a ${\bf F}_q$-rational divisor. Suppose that $$2g-2<\deg({\bf G})<n,$$  and $$supp({\bf G}) \bigcap {\bf P}=\emptyset.$$ Then an algebraic geometric functional code is defined by evaluating rational functions in the space $L({\bf G})$ at points $P_1,\ldots, P_n$, $${\bf C}({\bf X}, {\bf G}, {\bf P})=\{(f(P_1), \ldots,f(P_n)): f \in L({\bf G})\}.$$ Then the dimension is $$k=\deg({\bf G})-g+1,$$ and the minimum weight is at least $$d \geq n-\deg({\bf G}).$$

It is clear that $$({\bf C}({\bf X}, {\bf G}, {\bf P}))^{q-1} \subset {\bf C}({\bf X}, (q-1){\bf G}, {\bf P}),$$ then if $(q-1)\deg({\bf G})<\frac{n}{2}$, ${\bf C}({\bf X}, {\bf G}, {\bf P})$ is an intersecting code, see \cite{Borello}. We have the following result, from the Drinfeld-Vl\u{a}dut bound, see \cite{TV}.\\

{\bf Theorem 2.1.} {\em For any given square prime power $q^{2m}$, we have a linear $(q-1)$-power intersecting algebraic code ${\bf C}({\bf X}, {|bf G}, {\bf P})$ over ${\bf F}_{q^{2m}}$ if $\deg({\bf G}) < \frac{n}{2(q-1)}$. Then we have a sequence of $q-1$-power intersecting algebraic geometric $q^{2m}$-ary codes with the rate at least $\frac{1}{q-1} \cdot (\frac{1}{2}-\frac{1}{q^m-1})$.}\\

\section{Minimal codes via $(q-1)$-power intersecting codes}

In this section, we give the proof of Theorem 1.1.\\

We recall the concatenation of a linear $[N, K, D]_{q^m}$ code ${\bf C}_1$ over ${\bf F}_{q^m}$ and a linear $[n, m, d]_q$ code ${\bf C}_2$ over ${\bf F}_q$. The concatenated code is a linear $[Nn,Km, \geq Dd]_q$ code over ${\bf F}_q$. Let $\phi: {\bf F}_{q^m}={\bf F}_q^m \longrightarrow {\bf F}_q^n$ be a ${\bf F}_q$-linear injective mapping, such the image of $\phi$ is the code ${\bf C}_2$. For each codeword ${\bf c}=(c_0, c_1, \ldots, c_{N-1}) \in {\bf C}_1 \subset {\bf F}_{q^m}^N$, $c_i \in {\bf F}_{q^m}$, we associated a vector $\Phi({\bf c})=(\phi(c_0), \phi(c_1), \ldots, \phi(c_{N-1})) \in {\bf F}_q^{Nn}$. Then $\Phi: {\bf C}_1 \longrightarrow {\bf F}_q^{Nn}$ is a ${\bf F}_q$-linear injective mapping. The image the concatenated $[Nn, mK]_q$ code. Since there are at least $D$ indices $i$ such that $\phi(c_i)$ is nonzero, and each nonzero $\phi(c_i)$ has its Hamming weight $d$, the minimum weight of the concatenated code is at least $Dd$.\\

{\bf Theorem 3.1.} {\em Let ${\bf C}_1 \subset {\bf F}_{q^m}^N$ be a $q^m$-ary $K$-dimensional linear code such that ${\bf C}_1^{q-1}\subset {\bf F}_{q^m}^N$ is an intersecting code. We assume that there is an explicit $[n,m]_q$ minimal $q$-ary code ${\bf C}_2 \in {\bf F}_q^n$. Then we can construct a minimal linear $[Nn,Km]_q$ code in ${\bf F}_q^{Nn}$.}\\

{\bf Proof.} Let ${\bf c}_1=(c_{0,1}, \ldots, c_{N-1,1})$ and ${\bf c}_2=(c_{0,2}, \ldots, c_{N-1,2})$ be two nonzero codewords of the $q^m$-ary linear ${\bf C}_1$ such that, $supp(\Phi({\bf c}_1)) \subset supp(\Phi({\bf c}_2)).$ Then for each coordinate position $i$, $0 \leq i \leq N-1$, we have $supp(\phi(c_{i,1})) \subset supp(\phi(c_{i,2}))$. Since $\phi(c_{i,1})$ and $\phi(c_{i,2})$ are two nonzero codewords of the $q$-ary minimal ${\bf C}_2$, then there exists a nonzero element $\lambda_i \in {\bf F}_q^*$, such that $$c_{i,1}=\lambda_i c_{i,2}.$$

It is clear that $$\lambda_i^{q-1} c_{i,2}^{q-1}=c_{i,1}^{q-1}=c_{i,2}^{q-1},$$ since $\lambda_i \in {\bf F}_q^*$ and $\lambda_i^{q-1}=1$. Then we have two nonzero codewords ${\bf c}_1'$ and ${\bf c}_2'$ in ${\bf C}_1^{q-1}$, such that, $$supp({\bf c}_1) \subset supp({\bf c}_2),$$ and at the coordinate positions of $supp({\bf c}_1')$, the coordinates of these two nonzero codewords at this position are the same. Then ${\bf c}_2'-{\bf c}_1'$ is a nonzero codeword of ${\bf C}_1^{q-1}$ such that $supp({\bf c}_2'-{\bf c}_1') \bigcap supp({\bf c}_1')$ is empty. This is a contradiction to the assumption that ${\bf C}_1^{q-1}$ is intersecting. Then the concatenated code is minimal.\\

{\bf Corollary 3.1.} {\em Let ${\bf C}_1 \subset {\bf F}_{q^m}^N$ be a $q^m$-ary $K$-dimensional linear code such that ${\bf C}_1^{q-1}\subset {\bf F}_{q^m}^N$ is an intersecting code. We assume that there is an explicit $[n,m]_q$ minimal $q$-ary code ${\bf C}_2 \in {\bf F}_q^n$. Then we construct a strong blocking set in $P({\bf F}_q, Km-1)$ with the size $Nn$.}\\

{\bf Corollary 3.2.} {\em Let ${\bf C}_1$ be a $q^m$-ary rate $R_1$ linear code such that ${\bf C}_1^{q-1}$ is an intersecting code. We assume that there is an explicit rate $R_2$ minimal $q$-ary code ${\bf C}_2$. Then we construct a $q$-ary minimal code with the rate $R_1 R_2$.}\\

{\bf Corollary 3.3.} {\em Suppose that $q$ is an arbitrary prime power and $m \geq 2$ is a positive integer. We assume that there exists an explicit $[n,2m]_q$ minimal code. Then we construct a sequence of strong blocking sets $S_i \subset P({\bf F}_q, k_i-1)$, where $k_i$ goes to the infinity,  such that $$\lim \frac{|S_i|}{k_i}\leq \frac{n(q^m-1)(q-1)}{m(q^m-3)}.$$ Equivalently, we construct sequences of minimal $q$-ary codes with rate at least $\frac{m(q^m-3)}{n(q-1)(q^m-1)}$.}\\

{\bf Proof.} The conclusion follows from Theorem 2.1 and 3.1.\\

{\bf Corollary 3.4.} {\em For $q=2$, we construct a sequence of linear binary minimal codes with the rate $\frac{145}{961} \approx \frac{1}{6.6276}$ explicitly.}\\

{\bf Proof.} In the case of $q=2$, we use the $[31,10]_2$ minimal codes in \cite[Table 1]{Kurz} as inner code. Then the result follows from Theorem 3.1.\\

{\bf Corollary 3.5.} {\em For $q=3$, we construct a sequence of linear trifferent (ternary minimal) codes with the rate $\frac{1}{20}$ explicitly.}\\

{\bf Proof.} In the case of $q=3$, we use the $[15,4]_3$ minimal codes in \cite[Page 11]{Kurz} as inner code. Then the result follows from Theorem 3.1.\\

{\bf Corollary 3.6.} {\em For $q=4$, we construct a sequence of linear $4$-ary minimal codes with the rate at least $\frac{13}{405}$ explicitly.}\\

{\bf Proof.} In the case of $q=4$, we use the $[16,4]_4$ minimal codes in \cite[Table 1]{Kurz} as inner code. Then the result follows from Theorem 3.1.\\

{\bf Corollary 3.7.} {\em For $q=5$, we construct a sequence of linear $5$-ary minimal codes with the rate at least $\frac{11}{864}$ explicitly.}\\

{\bf Proof.} In the case of $q=5$, from the upper bound based on probabilistic argument, there is a linear $[36,4]_5$ minimal code. This code can be constructed explicitly, since the dimension is $4$. Using this code as an inner code, the result follows from Theorem 3.1.\\

\section{Short minimal codes over small fields}

We first give ssymptotic upper bounds for $\lim_{k \longrightarrow \infty} \frac{m(k,q)}{k}$ for $ q \leq 11$. As showed in Corollary 1.1, 1.2 and 1.3, our construction has smaller sizes for strong blocking sets over ${\bf F}_q$, for $q \leq 11$, when compared with \cite[Table]{Alon}. In the following table, we list our upper asymptotic bounds for $\lim_{k \longrightarrow \infty} \frac{m(k,q)}{k}$. The lower bound is from \cite[Table 1]{Scotti}.\\

\begin{longtable}{|l|l|l|l|l|}
\caption{\label{tab:6-1} Upper and lower bounds for $\lim_{k \longrightarrow \infty}\frac{m(k,q)}{k}$}\\ \hline
$q$     & Lower bound &Upper bound&Upper bound in \cite{Alon}&References \\ \hline
$2$           & $3.5276$ &$5.8334$&$354$ & \cite{Borello}\\ \hline
$3$       & $4.5516$ &$13.56$&$308$ & \cite{Bishnoi}\\ \hline
$4$    & $5.568$ &$31.1538$&$295$& This paper \\ \hline
$5$ &$6.5805$ &$78.5455$&$324$& This paper\\ \hline
$7$ &$8.5987$&$150.29$ & $376$&This paper\\ \hline
$8$&$9.6057$&$195.1967$ &$378$&This paper\\ \hline
$9$ &$no$& $246.1540$ &$390$&This paper\\ \hline
$11$ &$no$& $363.5593$ &$420$&This paper\\ \hline
\end{longtable}	

In the above table, we use the upper bound $$m(4,q) \leq \left\lceil \frac{2}{1+\frac{1}{(q+1)^2ln q}}3\right\rceil (q+1)$$ for $q \geq 3$. By a search algorithm, we can assume minimal codes of dimension $4$ satisfying this bound have been explicitly constructed. Then from our simple construction, the asymptotic bounds in the above table are obtained.\\

It is always interesting to construct shortest minimal of length $m(k,q)$ for given fixed $q$ and $k$. Only little is known for this problem, see \cite{Kurz,Kurz1,Scotti1}. In \cite[Page 6]{Scotti1}, $m(k,2)$, $2 \leq k \leq 9$ were listed, and $m(k,3)$ (trifferent code case), $2\leq k \leq 5$ were listed.  We have the following result for trifferent codes.\\

{\bf Proposition 4.1.} {\em Suppose that there exists an projective smooth irreducible curve defined over ${\bf F}_{81}$
of genus $g$. Let $N$ be the number of its ${\bf F}_{81}$-rational points. Then we construct an explicit $[14(N-t),4(u-g+1)]_3$ trifferent code for $u < \frac{N-t}{2}.$}\\

For example, there is an elliptic curve defined over ${\bf F}_{81}$, with $100$ rational points. Then we construct $[14(100-t), 4u]_3$ trifferent codes for $t=1,2,3,\ldots,$ and $u < \frac{100-t}{2}.$ That is, linear trifferent $[1386,196]_3$ code and linear trifferent $[1358,192]_3$ code are constructed explicitly.\\

\end{document}